\documentclass{article}
\usepackage{amssymb}
\newcommand{\qed}{\hfill $\Box$}
\newtheorem{teorema}{Theorem}[section]
\newtheorem{posledica}{Corollary}[section]

\newtheorem{lema}{Lemma}[section]

\title{Some Properties of Posynomial Rings}
\author{\v Zarko Mijajlovi\' c, Milo\v s Milo\v sevi\' c, Aleksandar Perovi\' c}
\date{}
\begin{document}
\maketitle
\begin{abstract}
In this article we shall study some basic properties of posynomial rings with particular emphasis on  rings
${\rm Pos}(K,\mathbb{Q})[\bar x]$, and ${\rm Pos}(K,\mathbb{Z})[\bar x]$.
The latter ring is the well known ring of Laurent polynomials.
\end{abstract}
\section{Introduction}

The notion of a posynomial\footnotetext{{\it Mathematical Subject Classification}: 16S34, 13A99

\hspace{1ex} {\it Key words}: group rings, Laurent polynomial rings, general commutative ring theory}
(positive polynomial) appeared in
geometric programming as a generalization of a polynomial. Zener
introduced posynomial functions about forty years ago in order to
compute minimal costs (see [7]). Aside from economy and
management, in the last decade posynomials have been used in
optimal integral circuit design (see [5], [6] and [8]).

The applicability of posynomials essentially relies on definability of root functions in the
theory of real closed fields (RCF) and on realtime procedures for quantifier elimination in RCF based
on the partial cylindrical algebraic decomposition.

We shall study here some algebraic and computational properties of rings of posynomials over a
commutative domain. In particular, it is proved that a posynomial ring ${\rm Pos}({\bf R},
{\bf G})[\bar x]$ is not noetherian and it is not UFD (unique factorization domain)
if ${\bf R}$ is a domain and ${\bf G}$ is an abelian group such that
$\bigcup_{n>1}G_n\not=\{0\}$, where $G_n=\bigcap_{k\in\mathbb{N}}n^kG$. Further, we introduce the posynomial
Zariski topology and prove the analogues to the Hilbert's Nullstellensatz and the real Nullstellensatz.
Finally, we shall study the ideal membership problem in the posynomial rings ${\rm Pos}({\bf K},\mathbb{Z})
[\bar x]$ and ${\rm Pos}({\bf K},\mathbb{Q})[\bar x]$ under assumption that $\bf K$ is a
computable domain.
\section{Preliminaries and notation}
Symbols $\mathbb{N},\mathbb{Z},\mathbb{Q},\mathbb{R}$ and $\mathbb{C}$ denote respectively the sets of
natural, integer, rational, real and complex numbers. Throughout this paper, we
assume that ${\bf R}=(R,+,\cdot,0,1)$ is a commutative domain with the multiplicative unit 1,
${\bf S}=(S,+,0)$ is a commutative semigroup, ${\bf G}=(G,+,0)$ is an abelian group and
${\bf K}=(K,+,\cdot,0,1)$ is a field.

For the given function $f:S\longrightarrow R$ we define its support by
$$
{\rm supp}(f)=\{s \in S\ |\ f(s)\not=0\}.
$$
The set of all functions $f:S\longrightarrow R$ with finite
supports we denote by $R[S]$.
If $f,g\in R[S]$ and $s\in S$, an addition and a multiplication on $R[S]$ are defined by
$$(f+g)(s)=f(s)+g(s),\ \ \  (fg)(s)=\sum_{u,v\in S,u+v=s}f(u)g(v).$$
If $\bf 0$ and $\bf 1$ are functions defined by
$$
{\bf 0}(s)=0,\ \ \ {\bf 1}(s)=\left\{
\begin{array}{lll}
1&,&s=0\\
0&,&s\not=0
\end{array}
\right. ,
$$
the structure ${\bf R}[{\bf S}]=(R[S],+,\cdot, {\bf 0},{\bf 1})$ is a commutative ring and it is
called a semigroup ring (see [2] and [13]).
\vspace{1ex}

The ideal $I$ of the ring $\bf R$ generated by $S\subseteq R$ will be denoted by $\langle
S\rangle_{\bf R}$; we omit $\bf R$ if the context is clear.
\vspace{1ex}

An ideal $I\subseteq R$ is real if
for each sequence $r_1,\dots,r_n$ of elements of $R$ we have that
if $r_1^2+\cdots +r_n^2\in I$ than each $r_i$ is in $I$. For the rest of notation
and definitions on real algebra we shall follow [3].
\vspace{1ex}

The dimension of $\bf R$ is the maximal length of strictly
increasing chains of prime ideals in $\bf R$.  More on dimension and integral
elements can be found in [9] and [11].
\section{Definition and basic properties}

We  introduce the notion of posynomial over $\bf R$ and $\bf S$ as a term of the form
\begin{eqnarray*}
\sum_{i=1}^n r_ix^{s_i},\ \ n\in\mathbb{N},r_i\in R,s_i\in S,
\end{eqnarray*}
where $x^0=1,\  x^{s_1}\cdot x^{s_2}=x^{s_1+s_2}$. The posynomial ring over $\bf R$ and $\bf S$ is
denoted by ${\rm Pos}({\bf R},{\bf S})[x]$, and we see that this ring is isomorphic to the semigroup ring
${\bf R}[{\bf S}]$. Posynomials in multiple variables are defined by induction:
$$
{\rm Pos}({\bf R},{\bf S})[x_1,\dots,x_{n+1}]={\rm Pos}({\rm Pos}
({\bf R},{\bf S})[x_1,\dots,x_n],{\bf S})[x_{n+1}].
$$
The following lemma is an easy fact on semigroup rings.
\begin{lema}
Let $\bf R$ be a commutative ring, let $\bf S$ be a commutative semigroup and suppose that $\bf S$ has
a finite cyclic subgroup.
Then the ring ${\rm Pos}({\bf R},{\bf S})[\bar x]$ is not a domain.
\end{lema}

\noindent Therefore, if $\bf S$ is a finite group or if $\bf S$ has an element of finite order, then
${\rm Pos}({\bf R},{\bf S})[\bar x]$ is not a domain.
\vspace{2ex}

Let ${\bf S}=(S,+,<,0)$ be an ordered semigroup. We say that a posynomial
$$
f(x)=\sum_{i=1}^n r_ix^{s_i},\ r_i\not=0
$$
is in ordered form if $s_1<\cdots <s_n$. In particular, let ${\rm deg}(f)=s_n$ be a degree of
the posynomial $f$.

\begin{lema}
Let $\bf R$ be a domain and let $\bf S$ be an ordered semigroup. Then the ring
${\rm Pos}({\bf R},{\bf S})[\bar x]$ is a domain.
\end{lema}
{\bf Proof}. Observe that the product of two monomials with nonzero
coefficients is not 0.
Let $f=r_1x^{s_1}+\cdots +r_nx^{s_n}$, $g=r'_1x^{s_1'}+\cdots +r_m'x^{s_m'}$,
$n>1$ or $m>1$,
$r_i,r_j'\not=0$, $s_1<\cdots<s_n$ and $s_1'<\cdots<s_m'$. Then
$$
fg=r_1r'_1x^{s_1+s_1'}+a_nr_m'x^{s_n+s_m'}\not=0
$$
since $s_1+s_1'<s_n+s_m'$. We use induction to complete the claim.
\hfill$\Box$

\begin{posledica}
Let the ring $\bf R$ be a domain and let $\bf G$ be a torsion free abelian group. Then
${\rm Pos}({\bf R},{\bf G})[\bar x]$ is a domain.
\end{posledica}
{\bf Proof}. Using the Malcev's compactness theorem one can prove that each torsion free abelian group
can be ordered, so by the previous lemma the claim follows.

\hfill$\Box$
\vspace{1ex}

Therefore,
${\rm Pos}({\bf R},{\bf G})[\bar x]$ is a domain if and only if the abelian group $\bf G$
is torsion free.

We use the same argument as in lemma 3.2 to prove:
\begin{teorema}
Let $\bf R$ be a domain and let $\bf G$ be an ordered abelian group. Then units in
${\rm Pos}({\bf R},{\bf G})[\bar x]$ are exactly monomials
$rx_1^{s_1}\cdots x_n^{s_n}$, where $r$ is an invertible element of $\bf R$.
\end{teorema}

For the given abelian group $\bf G$ and an integer $n>1$ let ${\bf G}_n=(G_n,+,0)$ be a subgroup of $\bf G$
defined by
$$
G_n=\bigcap_{k\in\mathbb{N}}n^kG.
$$
\begin{teorema}
Let $R$ be a domain and let $G$ be an ordered abelian group. If
$$
\bigcup_{n>1}G_n\not=\{0\},
$$
then ${\rm Pos}({\bf R},{\bf G})[\bar x]$ is not noetherian.
\end{teorema}
{\bf Proof}. Let $s\in \bigcup_{n>1}G_n\setminus \{0\}$. Then there are an integer $n>1$ and a sequence
$s_0, s_1, s_2,\dots$ in $G$ such that
$$
s=s_0=ns_1=n^2s_2=n^3s_3=\cdots.
$$
We claim that the chain
$$
\langle x_i^{s_0}-1\rangle\subseteq \langle x_i^{s_1}-1\rangle\subseteq \langle x_i^{s_2}-1\rangle
\subseteq\cdots
$$
is strictly increasing. Note that
$$
x_i^{s_n}-1=x_i^{ns_{n+1}}-1=(x_i^{s_{n+1}}-1)(x_i^{(n-1)s_{n+1}}+\cdots+1).
$$
Otherwise, let
$$
x_i^{s_{n+1}}-1=(x_i^{s_n}-1)\cdot f,
$$
where $f\in {\rm Pos}({\bf R},{\bf G})[\bar x]$. Then,
$$
x_i^{s_{n+1}}-1=(x_i^{s_{n+1}}-1)\cdot(x_i^{(n-1)s_{n+1}}+\cdots+1)\cdot f,
$$
which yields that
$$
(x_i^{(n-1)s_{n+1}}+\cdots+1)\cdot f=1.
$$
This is a contradiction, since $x_i^{(n-1)s_{n+1}}+\cdots+1$ is not a unit in the ring
${\rm Pos}({\bf R},{\bf G})[\bar x]$.
\hfill$\Box$
\vspace{2ex}

Note that converse implication doesn't hold. For instance, let $\bf G$ be a countable direct
sum of copies of $\mathbb{Z}$. Then
$$
\bigcup_{n>1}G_n=\{0\},
$$
since for each $s\in \mathbb{Z}$ we have that $|s|<n^{|s|}$.
${\rm Pos}({\bf R},{\bf G})[\bar x]$ is isomorphic to the ring of Laurent polynomials with
$\aleph_0$ variables, so it is not noetherian.
\vspace{2ex}

By the proof of the previous theorem we can conclude that ${\rm Pos}({\bf R},\mathbb{Q})
[\bar x]$ does not satisfy the ACC for principal ideals, so it cannot be UFD nor noetherian.
\vspace{1ex}

Let
$f=\sum_{i=1}^k c_ix_1^{s_{i1}}\cdots x_n^{s_{in}}\in {\rm Pos}({\bf R},\mathbb{Z})[\bar x]$.
We define the polynomial $F(f)\in {\bf R}[\bar x]$ by
$$
F(f)=x_1^{\alpha_1}\cdots x_n^{\alpha_n}\cdot f,
$$
where $\alpha_i={\rm max}\{-s_{1i},\dots,-s_{ki}\}$.

Note that $F$ is compatible with $\cdot$ (i.e.
$F(fg)=F(f)F(g)$), but it is not compatible with $+$ (for instance $F(x+1)=x+1$ and $F(1)=F(x)=1$).
It is easy to see that $F(f)$ is irreducible in ${\bf R}[\bar x]$ if and only if  $f$ is irreducible
in ${\rm Pos}({\bf R},\mathbb{Z})[\bar{x}]$.
\vspace{1ex}

For an arbitrary positive integer $m$ let us define a ring monomorphism
$\Phi_m:{\rm Pos}({\bf R},\mathbb{Q})[\bar x]\longrightarrow {\rm Pos}({\bf R},
\mathbb{Q})[\bar x]$
by
$$
\Phi_m(\sum_{i=1}^k c_ix_1^{s_{i1}}\cdots x_n^{s_{in}})= \sum_{i=1}^k c_ix_1^{ms_{i1}}
\cdots x_n^{ms_{in}}.
$$
Further, if $f_1,\dots,f_k$ are arbitrary posynomials from ${\rm Pos}({\bf R},\mathbb{Q})[\bar x]$,
then let $\pi(f_1,\dots,f_k)$ be the least positive integer $m$ such that each $\Phi_m(f_i)$
is a Laurent polynomial. It is easy to see that
$\Phi_m(f)\in {\rm Pos}({\bf R},\mathbb{Z})[\bar x]$ iff
$\pi(f)|m$, and thus
$$
\pi(f_1,\dots,f_k)={\rm LCM}(\pi(f_1),\dots,\pi(f_k)).
$$
Let $f\in {\rm Pos}({\bf R},\mathbb{Q})[\bar x]$ and let $m=\pi(f)$. Then $f$ is atomic iff
for each positive integer $n$ the polynomial $F(\Phi_{mn}(f))$ is irreducible in ${\bf R}[\bar x]$.

For example, there are no atomic elements in ${\rm Pos}(\mathbb{R},\mathbb{Q})[x]$ and
${\rm Pos}(\mathbb{C},\mathbb{Q})[x]$, since each polynomial of degree greater than $2$ is
reducible in $\mathbb{R}[x]$, and each polynomial of degree greater than $1$ is not atomic
in $\mathbb{C}[x]$. On the other hand, the posynomial $x+2$ is atomic in ${\rm Pos}(\mathbb{Q},\mathbb{Q})[x]$,
since each polynomial $F(\Phi_n(x+2))=x^n+2$ is by Eisenstein criterion irreducible
in $\mathbb{Q}[x]$.

Since $\langle f_1,f_2\rangle_{{\rm Pos}({\bf R},\mathbb{Q})[x]}=
\langle g\rangle_{{\rm Pos}({\bf R},\mathbb{Q})[x]}$, where
$$F(\Phi_{\pi(f_1,f_2)}(g))={\rm GCD}(F(\Phi_{\pi(f_1,f_2)}(f_1)),F(\Phi_{\pi(f_1,f_2)}(f_2)),$$
we see that every finitely generated ideal in ${\rm Pos}({\bf K},\mathbb{Q})[x]$ can be
generated by one element.
\vspace{1ex}

\noindent {\bf Example.} The ideal
$I=\langle x^{\frac{1}{n}}-1\ |\ n\in \mathbb{N}\rangle_{{\rm Pos}({\bf R},\mathbb{Q})[x]}$
is prime:
\vspace{1ex}

suppose that $fg\in I$; then there is
a positive integer $n$ such that $fg\in \langle x^{\frac{1}{n}}-1\rangle_
{{\rm Pos}({\bf R},\mathbb{Q})[x]}$. Further, there is $h\in {\rm Pos}({\bf R},\mathbb{Q})[x]$
such that $fg=h(x^{\frac{1}{n}}-1)$. Let $m=\pi(f,g,h,x^{\frac{1}{n}}-1)$. Then
$$
F(\Phi_m(f))F(\Phi_m(g))=F(\Phi_m(h))F(\Phi_m(x^{\frac{1}{n}}-1)),
$$
so $x-1$ divides at least one of polynomials
$F(\Phi_m(f))$ and $F(\Phi_m(g))$; say $F(\Phi_m(f))$. We conclude that
$f\in\langle x^{\frac{1}{m}}-1 \rangle_{{\rm Pos}({\bf R},\mathbb{Q})[x]}$.

\begin{teorema}
If $\bf K$ is a field, then
$
{\rm dim} ({\rm Pos}({\bf K},\mathbb{Q})[x_1,\dots,x_n])=n .
$
\end{teorema}
{\bf Proof}. Note that for a nonzero integer $n$ each posynomial
$x_i^{\frac{1}{n}}$ is a zero of a monic polynomial
$y^{{\rm sgn}(n)n}-x_i^{{\rm sgn}(n)}$ over ${\rm Pos}({\bf K},\mathbb{Z})[\bar x]$,
 so ${\rm Pos}({\bf K},\mathbb{Q})[\bar x]$ is an integral
extension of ${\rm Pos}({\bf K},\mathbb{Z})[\bar x]$. Hence, the dimension of the
posynomial ring
${\rm Pos}({\bf K},\mathbb{Q})[\bar x]$
is equal to the dimension of ${\rm Pos}({\bf K},\mathbb{Z})[\bar x]$ and since dimension is a local property
we have that
$$
{\rm dim} ({\rm Pos}({\bf K},\mathbb{Q})[x_1,\dots,x_n])={\rm dim}({\bf K}[x_1,\dots,x_n])=n .
$$\hfill $\Box$

We observe that posynomials from ${\rm Pos}({\bf K},\mathbb{R})[\bar x]$ which annul some
polynomial with coefficients from ${\bf K}[\bar x]$ are exactly the elements of the ring
${\rm Pos}({\bf K},\mathbb{Q})[\bar x]$.
\vspace{2ex}

At the end of this section we discuss the possibility of functional representation of
posynomials with positive rational exponents. Let $\bf K$ be a finite field of prime characteristic
$p$. The inverse of the Frobenius automorphism $x\mapsto x^{p^n}$ is a unique function on $\bf K$ which
satisfies natural equalities for the $p^n$-th root function $\phi$:
$$
(\phi(x))^{p^n}=x\ \ \ \ \mbox{   and   }\ \ \ \ \phi(xy)=\phi(x)\phi(y).
$$
Let $\bf K$ be an algebraic extension of the prime field $\mathbb{Z}_p$. Since each $a\in K$ is
contained in some finite field $\bf L$, again we conclude that there is a unique $b\in K$
such that $b^{p^n}=a$ and the corresponding $p^n$-th root function is compatible with multiplication.
Note that the same is true for an algebraically closed field of characteristic $p$, since
the polynomial $x^{p^n}-a$ has exactly one zero in that field.

Thus, if $\bf K$ is an algebraically closed field of characteristic $p>0$ or an algebraic extension
of prime field $\mathbb{Z}_p$, then
each posynomial $f$ in one variable over $\bf K$ of the form
$$
a_1x^{\frac{l_1}{p^{n_1}}}+\cdots+ a_kx^{\frac{l_k}{p^{n_k}}},\ l_i,n_j\in \mathbb{N}
$$
has a natural functional representation ${\bf f}:K\longrightarrow K$. Further, if
$$
S=\{ \frac{l}{p^n}\ |\ l,n\in\mathbb{N}\},
$$
then by $f\mapsto \bf f$ is defined a ring homomorphism from ${\rm Pos}({\bf K},
{\bf S})[x]$ into $K^K$.
\vspace{2ex}

Observe that the functional representation of $x^{\frac{1}{n}}\in {\rm Pos}(\mathbb{C},
\mathbb{Q})[x]$ determined by some branch of the $n$-th root is not compatible with multiplication
in $\mathbb{C}$.
\section{Laurent polynomials}

Let $\bf K$ be an arbitrary field of characteristic 0. The ring of Laurent polynomials over $\bf K$
(in variables $x_1,\dots,x_n$)
is the ring ${\rm Pos}({\bf K},\mathbb{Z})[\bar x]$. Note that ${\rm Pos}
({\bf K},\mathbb{Z})[\bar x]$ is just the
localization of ${\bf K}[\bar x]$ at $x_1\cdots x_n$, so it is noetherian, UFD and a graded ring.
\vspace{2ex}

We define the Zariski topology for Laurent polynomials in a similar way as in the case of polynomial
Zariski topology. Let
$$
K^n_{\not=0}=\{ (a_1,\dots,a_n)\in K^n\ |\ a_1\cdots
a_n\not=0\}.
$$
Each Laurent polynomial $f=\sum_{i=1}^kc_1x_1^{s_{i1}}\cdots x_n^{s_{in}}$ defines an unique function
${\bf f}:K^n_{\not=0}\longrightarrow K$ in a quite natural way:
$$
{\bf f}(a_1,\dots,a_n)=\sum_{i=1}^kc_1a_1^{s_{i1}}\cdots a_n^{s_{in}}.
$$
Note also that the mapping $f\mapsto \bf f$ is an embedding of the ring ${\rm Pos}(
{\bf K},\mathbb{Z})[\bar x]$
into the ring $K^{K^n_{\not=0}}$.

Let $S\subseteq {\rm Pos}({\bf K},\mathbb{Z})[\bar x]$ be an arbitrary set of Laurent polynomials.
A posynomial set in
$K^n_{\not=0}$ generated by $S$ is the set
$$
V_{\rm Pos}(S)=\{
(a_1,\dots,a_n)\in K^n_{\not=0}\ |\ (\forall f\in S)f[\bar{a}]=0\}.
$$
First, let us observe that $K^n_{\not=0}$ is the Zariski open set in affine space $K^n$ given as the
complement of the Zariski closed set $V(x_1\cdots x_n)$. Further,
$$
V_{\rm Pos}(S)=\bigcap_{f\in S}V(F(f))\ \cap\ K^n_{\not=0},
$$
so the posynomial sets (which are the base closed sets in the
posynomial Zariski topology) are closed in the induced topology on
the open subset $K^n_{\not=0}$ of the Zariski topology on $K^n$.
Thus we can immediately conclude that $K^n_{\not=0}$ is a Frechet
space (in the posynomial Zariski topology) and each posynomial
function is continuous. Further, since each two nonempty Zariski
open sets meet each other, the same will obviously hold for each
two nonempty posynomial Zariski open sets, thus $K^n_{\not=0}$ is
not a Hausdorff space. The compactness of $K^n_{\not=0}$ can be shown exactly in
the same way as for $K^n$ with polynomial Zariski topology.
\vspace{3ex}

As dual notion to posynomial sets, for an arbitrary set $X\subseteq K^n_{\not=0}$ let
$$
I_{\rm Pos}(X)=\{ f\in {\rm Pos}(K,\mathbb{Z})[\bar{x}]\ |\ (\forall (a_1,\dots,a_n)\in X)f[\bar{a}]=0\}.
$$
The ring ${\rm Pos}({\bf K},\mathbb{Z})[\bar x]/_{I_{\rm Pos}(X)}$
is reduced.
In particular, $I_{\rm Pos}(X)$ is a radical ideal. The next two results are analogues
of the corresponding polynomial theorems. The argument is similar, so we give only the
proof of real Nullstellensatz.

\begin{teorema}[Nullstellensatz for Laurent Polynomials]
Let $\bf K$ be an algebraically closed field and let $I$ be an arbitrary ideal in
${\rm Pos}({\bf K},\mathbb{Z})[\bar x]$. Then $V_{\rm pos}(I)\not=\emptyset$
if and only if $I$ is a proper ideal.
\end{teorema}
\noindent {\bf Remark.} The Hilbert's Nullstellensatz does not hold in
${\rm Pos}(\mathbb{C},\mathbb{Q})[x]$.
\vspace{1ex}

First let us observe that the function which maps $x^{\frac{1}{n}}$ to the principal
branch of the $n$-th root function is a ring embedding of ${\rm Pos}(\mathbb{C},
\mathbb{Q})[x]$ into $\mathbb{C}^{\mathbb{C}_{\not=0}}$.

Then
$I(V(\langle x^{\frac {1}{2}}+1\rangle ))=I(\emptyset)={\rm Pos}(\mathbb{C},\mathbb{Q})[x]$,
but $1\notin {\rm rad}\langle x^{\frac {1}{2}}+1\rangle $.
\vspace{3ex}

\begin{teorema}[Real Nullstellensatz for Laurent polynomials]
Let $K$ be a real closed field and let $I$ be an ideal in ${\rm Pos}({\bf K},\mathbb{Z})[\bar x]$.
Then
$$
I=I_{\rm Pos}(V_{\rm Pos}(I))
$$
if and only if $I$ is a real ideal.
\end{teorema}
{\bf Proof}. We will consider only nontrivial direction.
Suppose that $I$ is a real ideal; then it is a radical ideal and can be represented as a finite intersection of prime ideals
$I_1,\dots,I_k$ in ${\rm Pos}({\bf K},\mathbb{Z})[\bar x]$.
Clearly,
$$
I\subseteq I_{\rm Pos}(V_{\rm Pos}(I)).
$$
Let $f\in I_{\rm Pos}(V_{\rm Pos}(I))\setminus I$;
for instance, let $f\notin I_1$. The ring ${\rm Pos}({\bf K},\mathbb{Z})[\bar x]$ is noetherian, so there
are $f_1,\dots f_k\in I_1$ such that $I_1=\langle f_1,\dots,f_k\rangle $.
Since each prime ideal is real, the field
$$
{\bf K}_1=\mathbb{Q}({\rm Pos}({\bf K},\mathbb{Z})[\bar x]/I_1)
$$
is real. Let ${\bf K}_2$ be a real closure of ${\bf K}_1$.
Each $x_i$ is invertible in
${\rm Pos}({\bf K},\mathbb{Z})[\bar x]$, so $x_i+I_1\not=I_1$ and $(x_1+I_1,\dots, x_n+I_1)$
is a witness for
$$
{\bf K}_2\models \exists \bar v(F(f)(\bar v)\not=0\ \land\
\bigwedge_{i=1}^n F(f_i)(\bar v)=0\ \land\ \bigwedge_{i=1}^n v_i\not=0).
$$
The submodel completeness of the theory of
real closed fields yields
$$
{\bf K}\models \exists \bar v(F(f)(\bar v)\not=0\ \land\
\bigwedge_{i=1}^n F(f_i)(\bar v)=0\ \land\ \bigwedge_{i=1}^n v_i\not=0),
$$
which contradicts the fact that $V_{\rm Pos}(f)\supseteq V_{\rm Pos}(I)
\supseteq V_{\rm Pos}(I_1)$.
\hfill$\Box$
\section{Posynomials over computable fields}

From now on we will assume that $\bf K$ is a computable field of characteristic 0.
\begin{lema}
Let $p_0,p_1,\dots, p_n$ be arbitrary distinct prime numbers and let $f_i=x^{\frac{1}{p_i}}-1$.
Then
$$
f_0\notin \langle f_1,\dots,f_n\rangle_{{\rm Pos}({\bf K},\mathbb{Q})[x]}.
$$
\end{lema}
{\bf Proof}. Otherwise, there are posynomials $g_1,\dots,g_k\in {\rm Pos}({\bf K},\mathbb{Q})[x]$ such
that
$$
f_0=g_1f_1+\cdots+g_nf_n.
$$
Let $m=\pi(f_0,f_1,\dots,f_n,g_1,\dots,g_n)$. Then there are unique positive integers $s$ and $d$ such that
$m=p_0^sd$ and ${\rm GCD}(p_0,d)=1$. For an arbitrary $i>0$ we have that
$$
\Phi_m(f_i)=x^{\frac{p_0^sd}{p_i}}-1=(x^{p_0^s}-1)(x^{p_0^s(\frac{d}{p_i}-1)}+x^{p_0^s(\frac{d}{p_i}-2)}
+\cdots+1)
$$
and each $\frac{d}{p_i}-j$ is an integer, so $\Phi_m(f_i)$ is divisible by $x^{p_0^s}-1$ in
${\rm Pos}({\bf K},\mathbb{Z})[x]$. But $\Phi_m(f_0)$ is not divisible by $x^{p_0^s}-1$,
and we obtain a contradiction. \hfill $\Box$
\vspace{1ex}

We see that $x^{1\over p_0}-1$ is not a member of the posynomial ideal generated by the set
$$
B=\{x^{1\over p}-1\ |\ p\in A \},
$$
where $p_0\notin A$ and each member of $A$ is a prime number. This is a consequence of the fact that
for each ideal $I$, $a\in I$ if and only if $a$ can be represented as a finite sum of the form
$\sum_{i=1}^k b_ia_i$, where $a_i$ belong to the set of generators for $I$.

\begin{teorema}
The problem of ideal membership in ${\rm Pos}({\bf K},\mathbb{Q})[x]$
(for the given computable field $\bf K$) is not decidable,
i.e. there is a nonrecursive ideal in the ring ${\rm Pos}({\bf K},\mathbb{Q})[x]$.
\end{teorema}
{\bf Proof}. Let $A$ be a nonrecursive subset of $\mathbb{N}$ and let $I$ be a posynomial ideal generated
by the set
$$
B=\{x^{1\over p_i}-1\ |\ i\in A\},
$$
where $p_0,p_1,p_2,\dots$ is an increasing enumeration of prime numbers. Then, by the previous lemma
$$
x^{1\over p_i}-1\in I\ \ {\rm if\ and\ only\ if}\ \ i\in A.
$$
So, any algorithm which decides the predicate $``x^{1\over p_i}-1\in I"$ will also decide the
predicate $``i\in A"$ contradicting the fact that $A$ is a nonrecursive set. \hfill \qed
\vspace{2ex}

In the rest of this section we will describe one test for the  membership
to finitely generated ideals in ${\rm Pos}({\bf K},\mathbb{Q})
[\bar x]$.
\begin{teorema} Let $\bf K$ be a computable field. The question of ideal membership in the ring of Laurent
polynomials ${\rm Pos}({\bf K},\mathbb{Z})[\bar x]$ is decidable. Moreover, there is an algorithm
for testing the membership to finitely generated ideals in ${\rm Pos}({\bf K},\mathbb{Q})
[\bar x]$.
\end{teorema}
{\bf Proof}.
Let $I=\langle f_1,\dots,f_n\rangle_{{\rm Pos}({\bf K},\mathbb{Z})[\bar x]}$
be an ideal in Pos$({\bf K},\mathbb{Z})[\bar x]$.
We notice that:
$$
\begin{array}{ll}
g\in \langle f_1,\ldots ,f_k\rangle_{{\rm Pos}({\bf K},\mathbb{Z})[\bar x]}\ {\rm if\ and\ only\ if}& \\
& (*)\\
(x_1\cdots
x_n)^\lambda F(g)\in \langle F(f_1),\ldots ,F(f_k)\rangle_{{\bf K}[\bar
x]},\ {\rm for\ some}\ \lambda \in \mathbb{N}.
\end{array}
$$
We can write (*) using the saturation ideal of
$\langle F(f_1),\dots, F(f_n)\rangle_{{\bf K}[\bar x]}$ by $x_1\cdots x_n$:
$$
g\in \langle f_1,\ldots ,f_k\rangle_{{\rm Pos}({\bf K},\mathbb{Z})[\bar x]}\
\Leftrightarrow\ F(g)\in \langle F(f_1),\ldots
,F(f_k)\rangle_{{\bf K}[\bar x]}\ :\ (x_1\cdots x_n)^\infty
$$
If $J$ is an ideal in ${\bf K}[x_1,\dots, x_n]$, $h\in {\bf K}[\bar x]$ and $y$ a new variable, then
$$
J:h^\infty=\langle J,1-yh \rangle\cap {\bf K}[\bar x],
$$
where $J:h^\infty$ is an ideal in ${\bf K}[\bar x]$ and $\langle J,1-yh \rangle$ is an ideal in
${\bf K}[\bar x,y]$. The Gr\" obner basis of $J:h^\infty$ with respect to the lexicographical order
$x_1<\cdots<x_n$ is equal to the intersection of ${\bf K}[\bar x]$ and the Gr\" obner basis
$B$ of $\langle J,1-yh \rangle$ with respect to the lexicographical order $x_1<\cdots<x_n<y$
(see [2]).
\vspace{1ex}

Now we have an algorithm for testing whether
$g\in \langle f_1,\ldots ,f_k\rangle_{{\rm Pos}({\bf K},\mathbb{Z})[\bar x]}$ or not:

First we will find the Gr\"{o}bner basis
$B$ (with the respect to the lexicographical order)
of $\langle F(f_1),\ldots ,F(f_n),1-yx_1\cdots x_n\rangle \subseteq K[\bar x,y];\ B\cap
K[\bar x]=B_1$ will be the Gr\"{o}bner basis of
$$
\langle F(f_1),\ldots ,F(f_n)\rangle \ :\ (x_1\cdots x_n)^\infty.
$$
We divide $F(g)$ by $B_1$ in the lexicographical order; if the
remainder is $0$ then $g\in \langle f_1,\ldots ,f_k\rangle_{{\rm
Pos}({\bf K},\mathbb{Z})[\bar x]}$, otherwise $g\notin \langle f_1,\ldots
,f_k\rangle_{{\rm Pos}({\bf K},\mathbb{Z})[\bar x]}$.
\vspace{3ex}

We prove the existence of a procedure for testing ideal membership to
finitely generated ideal in ${\rm Pos}({\bf K},\mathbb{Q})[\bar{x}]$, where $\bf K$ is a computable
field. Let the ideal $J\subseteq {\rm
Pos}({\bf K},\mathbb{Q})[x_1,\ldots,x_n]$ be generated by
$f_1,\ldots,f_m$ and let $g$ be an arbitrary posynomial in ${\rm
Pos}({\bf K},\mathbb{Q})[\bar{x}]$.

Suppose that $g\in J$. There exist $h_1,\ldots,h_m\in {\rm
Pos}({\bf K},\mathbb{Q})[\bar{x}]$ such that
\begin{equation}
g=h_1f_1+\ldots+h_mf_m.
\end{equation}
We will write down all exponents which appear in
$g,f_1,\ldots,f_m$ in the form $\frac{p_i}{q_i},\ {\rm GCD}(p_i,q_i)=1,\
p_i\in \mathbb{Z},\ q_i\in \mathbb{N}$, the exponents which appear
in $h_1,\ldots,h_m$ in the form $\frac{k_i}{l_i}$, where $k_i$ and
$l_i$ are relatively prime,  and we will denote the least common
multiple of denominators $q_i$ by $s$ (note that $s=\pi(g,f_1,\dots,f_m)$). Now, we rewrite exponents
$\frac{p_i}{q_i}$ in the form $\frac{t_i}{s}$. Assume that the
posynomial $h_1$ contains a monomial $M$ with variable $x_i$
to the power $\frac{k_{i_0}}{l_{i_0}},\ l_{i_0}\nmid s.$ Then, the
product $h_1f_1$ contains monomial $M_1$ with variable $x_i$
to the power
$$a=\frac{k_{i_0}}{l_{i_0}} +\frac{t_{j_0}}{s}=\frac{k_{i_0}s+l_{i_0}t_{j_0}}{l_{i_0}s}.$$
Since $l_{i_0}\nmid s$ and ${\rm GCD}(k_{i_0},l_{i_0})=1$ we conclude
that $a$ is not of the form $\frac{t}{s}$ and that the monomial
$M_1$ cannot appear on the left side of the equation (1). We thus
obtain that all monomials with the same property as $M_1$ must
cancel and that $g$ can be expressed as
$$g=\tilde{h}_1f_1+\ldots+\tilde{h}_mf_m,$$
where all denominators $l_i$ of exponents $\frac{k_i}{l_i}$ which
occur in $\tilde{h}_1, \ldots,\tilde{h}_m$ divide $s$.
Thus
$$g\in \langle f_1,\ldots,f_m\rangle_{{\rm Pos}({\bf K},\mathbb{Q})[\bar{x}]}
\mbox{  iff  } \Phi_s(g)\in \langle
\Phi_s(f_1),\dots,\Phi_s(f_m)\rangle_ {{\rm
Pos}({\bf K},\mathbb{Z})[\bar{x}]}.$$
\hfill $\Box$

\vspace{1ex}

\begin{flushleft}
\v ZARKO MIJAJLOVI\' C\\ FACULTY OF MATHEMATICS\\ UNIVERSITY OF BELGRADE\\
STUDENTSKI TRG 16, 11000 BEOGRAD\\ SERBIA AND MONTENEGRO \\ \vspace{1ex}
{\it E-mail:}\ {\bf zarkom@eunet.yu}\\
\vspace{1ex}
MILO\v S MILO\v SEVI\' C\\
MATHEMATICAL INSTITUTE\\ SERBIAN ACADEMY OF SCIENCES AND ARTS\\
KNEZA MIHAILA 35, 11001 BEOGRAD\\ SERBIA AND MONTENEGRO \\ \vspace{1ex}
{\it E-mail:}\ {\bf mionamil@eunet.yu}\\ \vspace{1ex}
ALEKSANDAR PEROVI\' C\\
MATHEMATICAL INSTITUTE\\ SERBIAN ACADEMY OF SCIENCES AND ARTS\\
KNEZA MIHAILA 35, 11001 BEOGRAD\\ SERBIA AND MONTENEGRO \\ \vspace{1ex}
{\it E-mail:}\ {\bf peramail314@yahoo.com}\\
\end{flushleft}

\end{document}